\documentclass[10pt,a4paper]{article}
\usepackage[a4paper, total={6in, 8in}]{geometry}
\usepackage[latin1]{inputenc}
\usepackage{amsmath}
\usepackage{amsfonts}
\usepackage{authblk}
\usepackage{amssymb}
\usepackage{graphicx}
\newcommand{\eqnref}[1]{(\ref{#1})}
\newcommand{\bx}{\boldsymbol{x}}
\newcommand{\bu}{\boldsymbol{u}}
\newcommand{\be}{\boldsymbol{e}}
\newcommand{\bzero}{\boldsymbol{0}}
\newcommand{\bone}{\boldsymbol{1}}
\author{Johannes Lutzeyer}
\author{Edward A. K. Cohen}
\affil{Imperial College London}
\date{}
\title{Correcting the estimator for the mean vectors in a multivariate errors-in-variables regression model}
\begin{document}
\maketitle
\begin{abstract}
The multivariate errors-in-variables regression model is applicable when both dependent and independent variables in a multivariate regression are subject to measurement errors. In such a scenario it is long established that the traditional least squares approach to estimating the model parameters is biased and inconsistent. The generalized least squares, ordinary least squares and maximum likelihood estimators (under the assumption of Gaussian errors) were derived in the seminal paper of \cite{gleser81}. However, the ordinary least squares and maximum likelihood estimators for the mean vectors were incorrectly derived. In this short paper we amend this error, presenting the correct estimators of the mean vectors. 
\end{abstract}

\section{Introduction}
\cite{gleser81} presented the multivariate errors-in-variables regression model. We recap the formulation here and borrow the notation. 

In a multivariate errors-in-variables model there are $n$ observed random vectors $\bx_i = (\bx_{1i}',\bx_{2i}')'$, where $\bx_{1i}$ is $p\times 1$ and $\bx_{2i}$ is $r\times 1$, $i=1,2,...,n$. The assumed model is
$$
\bx_i = \left(\begin{array}{c} \bx_{1i}
 \\ 
\bx_{2i}
\end{array} \right) = \left(\begin{array}{c} \bu_{1i}
 \\ 
\bu_{2i}
\end{array} \right) + \left(\begin{array}{c} \be_{1i}
 \\ 
\be_{2i}
\end{array} \right) = \bu_i + \be_i
$$
$$
\bu_{2i} = \alpha + B\bu_{1i}\hspace{2cm}i=1,2,...,n.
$$
Here, the $(p+r)$-dimensional random vectors $\be_i$, $i=1,2,...,n$ are i.i.d. with mean $\bzero$ and covariance matrix $\Sigma_e$. The vectors $\bu_{1i}$, $i=1,2,...,n$, and $\bu_{2i}$, $i=1,2,...,n$ are referred to as the {\it mean vectors}. The transformation parameters are the $r\times p$ matrix $B$ and $r$-dimensional vector $\alpha$. There are two distinct cases. The {\it no-intercept} model assumes $\alpha$ is known to be $\bzero$ and $B$ and $\bu_{1i}$, $i=1,2,...,n$ are unknown and to be estimated. The more general {\it intercept} model has $\alpha$, $B$ and $\bu_{1i}$, $i=1,2,...,n$ as unknown and to be estimated. 

As in \cite{gleser81}, notation is conveniently condensed by letting $X_1$ be the $p\times n$ matrix with columns $\bx_{1i}$, $i=1,2,...,n$, letting $X_2$ be the $r\times n$ matrix with columns $\bx_{2i}$, $i=1,2,...,n$, letting $U_1$ be the $p\times n$ matrix with columns $\bu_{1i}$, $i=1,2,...,n$, and letting $U_2$ be the $r\times n$ matrix with columns $\bu_{2i}$, $i=1,2,...,n$. We also let $E$ be the $(p+r)\times n$ random matrix whose columns are $\be_i$, $i=1,2,...,n$. By constructing $X = (X_1',X_2')'$ and $U = (U_1',U_2')'$, we can represent the model as
$$
X = U + E = \left(\begin{array}{c}\bzero
 \\ 
\alpha\bone'_n
\end{array} \right)
+\left(\begin{array}{c}I_p
 \\ 
B
\end{array} \right)U_1 + E,
$$
where $\bone_n$ is an $n$-dimensional column vector of ones. \cite{gleser81} makes the following assumption:
\begin{center}
the columns of $E$ are i.i.d. with common mean vector $\bzero$ and common covariance matrix $\Sigma_e = \sigma^{2}I_{p+r}$, where scalar $\sigma^2>0$ is unknown,
\end{center}
although Section 5 of \cite{gleser81} demonstrates how this can be relaxed to a more general assumption of $\Sigma_e = \sigma^2\Sigma_0$, $\Sigma_0$ known.

\cite{gleser81} presents three approaches to estimating the unknown parameters $U_1$, $\alpha$, $B$ and $\sigma^2$ of the multivariate ``errors-in-variables'' model. Firstly, in the {\it maximum likelihood} approach the columns of $E$ are assumed i.i.d. multivariate normal allowing estimation of all four parameters. This estimation procedure is referred to as MLE. Secondly, in the {\it ordinary least squares} approach $U_1$, $\alpha$ and $B$ are chosen to minimize an orthogonally invariant norm of the residual matrix
$$
R(\alpha,B,U_1;X) = X-\left(\begin{array}{c}
\bzero \\ \alpha\bone_n' 
\end{array} \right)-\left(\begin{array}{c}
I_p \\ B 
\end{array} \right)U_1.
$$
This estimation procedure is referred to as OLSE. Thirdly, in the {\it generalized least squares} approach $\alpha$ and $B$ are chosen to minimize an orthogonally invariant norm of the normalized residual matrix
$$
Q(\alpha,B;X) = (I_r + BB')^{-1/2}(X_2 - \alpha\bone_n' - BX_1),
$$
where $(I_r + BB')^{1/2}$ is any square root of $I_r + BB'$. This estimation procedure is referred to as GLSE. 

The benefit of GLSE and OLSE approaches is they make no assumption on the distribution of the columns of $E$ beyond those given for their mean vector and covariance matrix.
\section{Correction}
In this section it becomes necessary to define the following terms in the same way as Section 2 of \cite{gleser81}. Let 
\begin{align*}
C_n & = I_n, \hspace{2.93cm} \text{for the no-intercept model}, \\
& = I_n - n^{-1}\bone_n\bone_n', \hspace{1cm} \text{for the intercept model},
\end{align*}
let
$$
W = XC_n X'
$$
and let $GDG'$ be the eigen-decomposition of $W$ such that $D$ is a diagonal matrix of the ordered eigenvalues and $(p+r)\times (p+r)$ matrix $G$ is orthogonal. We partition $G$ as 
$$
\left(\begin{array}{cc}
G_{11} & G_{12} \\ 
G_{21} & G_{22}
\end{array} \right),
$$
where $G_{11}$ is $p\times p$.

\cite{gleser81} presents the following estimators:
\begin{equation}
\label{B}
\hat{B} = G_{21}G_{11}^{-1},
\end{equation}
\begin{align}
\label{alpha}
\hat{\alpha} & = \bzero,\hspace{2.4cm}\text{for the no-intercept model,} \\ & = (-\hat{B},I_r)\bar{\bx},\hspace{1cm}\text{for the intercept model,}\nonumber
\end{align}
\begin{equation}
\hat{U}_1(\hat{\alpha},\hat{B}) = \left( G_{11} G_{11}' ~ X_1 + G_{11} G_{21}' ~ X_2  \right) C_n.\label{gleser_estim}
\end{equation}
where $\bar{\bx} = n^{-1}X\bone_n$. In Theorem 2.1 in \cite{gleser81} (\ref{B}) and (\ref{alpha}) are shown to be the OLSEs and GLSEs for $B$ and $\alpha$, and (\ref{gleser_estim}) the OLSE for the mean vectors. For the \underline{no-intercept} model this is correct. However, for the \underline{intercept} model, while $\hat{B}$ and $\hat{\alpha}$ are correct and it is correctly identified that the OLSE of the matrix of mean vectors $U_1$ is given by
\begin{equation}
\hat{U}_1(\hat{\alpha},\hat{B}) = \left(I_p + \hat{B}' \hat{B}\right)^{-1} \left( I_p, \hat{B}' \right) \left( X - 
\begin{pmatrix}
	\mathbf{0}\\
	\hat{\alpha}\mathbf{1}_n'
\end{pmatrix}
\right),
\label{gleser_min}
\end{equation}
it is then wrongly stated that this is equal to the key expression given in \eqnref{gleser_estim}.

In the following the term in \eqnref{gleser_min} is manipulated to show it does not lead to the estimator in \eqnref{gleser_estim}. Expanding \eqnref{gleser_min} we get
\begin{align}
\hat{U}_1(\hat{\alpha},\hat{B}) &= \left(I_p + \hat{B}' \hat{B}\right)^{-1} \left( X_1 + \hat{B}' (X_2 - n^{-1}(-\hat{B}X_1 + X_2)~ \mathbf{1}_n \mathbf{1}_n') \right)\nonumber\\
&=\left(I_p + \hat{B}' \hat{B}\right)^{-1}\times\nonumber\\ &\hspace{1cm}\left( X_1 + n^{-1}\hat{B}' \hat{B}~ X_1  ~\mathbf{1}_n \mathbf{1}_n' +
\hat{B}' X_2 -n^{-1}\hat{B}' X_2 ~\mathbf{1}_n \mathbf{1}_n') \right). \label{gleser_eqn1}
\end{align}
Since $G$ is orthogonal, considering the top left block of the equality $G' G = I_{p+r}$ we get the relation $G_{11}' G_{11} + G_{21}' G_{21}= I_p$. Using this and $\left(G_{11}^{-1}\right)' = \left(G_{11}'\right)^{-1}$ the following simplification arises:
\begin{align}
\hat{B}' \hat{B} = \left(G_{11}^{-1}\right)' G_{21}' G_{21} G_{11}^{-1} & =  \left(G_{11}^{-1}\right)' \left( I_p - G_{11}' G_{11}\right) G_{11}^{-1}\nonumber \\ & = \left(G_{11}^{-1}\right)' G_{11}^{-1} -I_p. \label{gleser_eqn2}
\end{align}
From \eqnref{gleser_eqn2} we have \eqnref{gleser_eqn1} simplify to:
\begin{multline*}
\hat{U}_1(\hat{\alpha},\hat{B}) = G_{11} G_{11}' \left( X_1 + n^{-1}\left( \left(G_{11}^{-1}\right)' G_{11}^{-1} -I_p \right) ~ X_1  \mathbf{1}_n \mathbf{1}_n'\right. \\ +\left.
\left(G_{11}'\right)^{-1} G_{21}' X_2 -n^{-1}\left(G_{11}'\right)^{-1} G_{21}' X_2 \mathbf{1}_n \mathbf{1}_n') \right),
\end{multline*}
giving the final result
\begin{equation}
\label{correct}
\hat{U}_1(\hat{\alpha},\hat{B}) = n^{-1} X_1\mathbf{1}_n \mathbf{1}_n' + \left( G_{11} G_{11}' ~ X_1 + G_{11} G_{21}' ~ X_2  \right) C_n.
\end{equation}
We are therefore left with the additional term $n^{-1} X_1 \mathbf{1}_n \mathbf{1}_n'$, which is not present in \eqnref{gleser_estim}. Using the correct expression (\ref{correct}), Theorem 2.3 in \cite{gleser81}, which states this to also be the MLE, is now also correct.

\bibliographystyle{apalike}
\bibliography{references}

\end{document}